\input amstex
\input Amstex-document.sty

\pageno 599

\def\shorttitle{Spectral Problems}
\def\shortauthor{A. Klyachko}

\topmatter %
\title\nofrills{\boldHuge Vector Bundles, Linear Representations, and Spectral Problems}
\endtitle

\author \Large Alexander Klyachko* \endauthor

\thanks *Department of Mathematics, Bilkent University,
Bilkent 06533, Ankara, Turkey. E-mail: klyachko\@fen.bilkent.edu.tr
\endthanks

\abstract\nofrills \centerline{\boldnormal Abstract}

\vskip 4.5mm

This paper is based on my talk at ICM on recent progress in a
number of classical problems of linear algebra and representation
theory, based on new approach, originated from geometry of stable
bundles and geometric invariant theory.

\vskip 4.5mm

\noindent {\bf 2000 Mathematics Subject Classification:} 14F05, 14M15, 14M17, 14M25,
15A42.

\noindent {\bf Keywords and Phrases:} Bundles, Linear representations,
Spectral problems.
\endabstract
\endtopmatter

\document

\baselineskip 4.5mm \parindent 8mm

\specialhead \noindent \boldLARGE 1. Introduction \endspecialhead
Theory of vector bundles brings a new meaning and adds a delicate
geometric flavour to classical spectral problems of linear
algebra, relating them to geometric invariant theory,
representation theory, Schubert calculus, quantum cohomology,  and
various moduli spaces. The talk may be considered as a supplement
to that  of Hermann Weyl \cite{35} from which I borrow the
following quotation

{\it ``In preparing this lecture, the speaker has assumed that he
is expected to talk on a subject in which he had some first-hand
experience through his own work. And glancing back over the years
he found that the one topic to which he has returned again and
again is the problem of eigenvalues and eigenfunctions in its
various ramifications."}

\specialhead \noindent \boldLARGE 2. Spectra and representations
\endspecialhead

Let's start with two classical and apparently independent
problems.

\proclaim{\bf\quad\quad Hermitian spectral problem} {\rm Find all
possible spectra $\lambda(A+B)$ of sum of Hermitian operators
$A,B$ with given spectra}
$$\aligned
\lambda(A):&\,\,\,\lambda_1(A)\ge\lambda_2(A)\ge\cdots\ge\lambda_n(A),\\
\lambda(B):&\,\,\,\lambda_1(B)\ge\lambda_2(B)\ge\cdots\ge\lambda_n(B).
\endaligned$$
\endproclaim
Among commonly known restrictions on spectra are {\it trace
identity}
$$\sum_i\lambda_i(A+B)=\sum_j\lambda_j(A)+\sum_k\lambda_k(B)$$
and a number of classical inequalities, like that of Weyl
\cite{34}
$$\lambda_{i+j-1}(A+B)\le \lambda_i(A)+\lambda_j(B).\tag2.0$$
%\roster
%\item"(ii)"

\proclaim{\quad\quad Tensor product problem}  {\rm Find all
components $V_\gamma\subset V_\alpha\otimes V_\beta$ of tensor
product of two irreducible representations of $\text{\rm GL}_n$
with highest weights (=Young diagrams)}
$$\aligned
\alpha:&\,\,\, a_1\ge a_2\ge\cdots\ge a_n\\
\beta:&\,\,\,  b_1\ge b_2\ge\cdots\ge b_n.
\endaligned
$$
\endproclaim
In contrast to the spectral problem (2.1) the coefficients of
tensor product decomposition
$$V_\alpha\otimes V_\beta=\sum_\gamma c_{\alpha\beta}^\gamma
V_\gamma\tag2.1$$ can be evaluated algorithmically by {\it
Littlewood--Richardson rule}, which may be described as follows.
Fill $i$-th row of diagram $\beta$ by symbol $i$. Then
$c_{\alpha\beta}^\gamma$ is equal to number of ways to produce
diagram $\gamma$ by adding cells from $\beta$ to $\alpha$ in such
a way that the symbols
\roster
\item"i)" weakly increase in rows,
\item"ii)" strictly increase in columns,
\item"iii)" reading all the symbols from right to left, and from
top to bottom produces a {\it lattice permutation}, i.e. in every
initial interval symbol $i$ appears at least as many times as
$i+1$.
\endroster

It turns out that these two problems are essentially equivalent
and have the same answer. To give  it, let's associate with a
subset $I\subset\{1,2,\ldots,n\}$ of cardinality $p=|I|$ Young
diagram $\sigma_I$ in a rectangular of format $p\times q$,
$p+q=n$, cut out by polygonal line $\Gamma_I$, connecting $SW$ and
$NE$ corners of the rectangular, with $i$-th unit edge running to
the North, for $i\in I$, and to the East otherwise. One can
formally multiply the diagrams by L-R rule
$$\sigma_I\sigma_J=\sum_k c^K_{IJ}\sigma_K\tag2.2$$
where $\sigma^K_{IJ}:=c^{\sigma_K}_{\sigma_I\sigma_J}$ are L-R
coefficients. Geometrically (2.2) is decomposition of product of
two Schubert cycles in cohomology ring of Grassmannian
$\text{G}^q_p$ of linear subspaces of dimension $p$ and
codimension $q$.

\proclaim{Theorem~2.1}  The following conditions are equivalent

{\rm i)} There exist Hermitian operators $A$, $B$, $C=A+B$ with
spectra $\lambda(A)$, $\lambda(B)$, $\lambda(C)$.

{\rm ii)} Inequality
$$\lambda_K(C)\le\lambda _I(A)+\lambda_J(B),
\tag{IJK}$$
holds each time L-R coefficient $c_{IJ}^K\ne0$. Here $I,J,K\subset
\{1,2,\ldots,n\}$ are subsets of the same cardinality $p<n$, and $
\lambda_I(A)=\sum_{i\in I}\lambda_i(A)$.

{\rm iii)} For integer spectra $\alpha=\lambda(A)$,
$\beta=\lambda(B)$, $\gamma=\lambda(C)$ the above conditions are
equivalent to
$$V_\gamma\subset V_\alpha\otimes V_\beta.\tag2.3$$
\endproclaim

\subhead{\bf Remarks~2.2}\endsubhead (1) The last claim iii) implies a {\it
recurrence procedure\,} to generate all $\alpha,\beta,\gamma$ with
$c_{\alpha\beta}^\gamma\ne 0$:
$$c_{\alpha\beta}^\gamma\ne
0\underset{\text{LR}}\to{\Longleftrightarrow} V_\gamma\subset
V_\alpha\otimes V_\beta
\underset{\text{Th}}\to{\Longleftrightarrow}
\gamma_K\le\alpha_I+\beta_J\text{ each time }c^K_{IJ}\ne0.$$
Here $c_{\alpha\beta}^\gamma$ are Littlewood-Richardson
coefficients for group $\text{GL}_n$, while $c^K_{IJ}$ are L-R
coefficient for group $\text{GL}_p$ of {\it smaller\,} rank $p<n$.
An explicit form of this recurrence has been conjectured by
A.~Horn \cite{13} in the framework of Hermitian spectral problem.

(2) Inequalities (IJK) for $c_{IJ}^K\ne 0$ define a cone in the
space of triplets of spectra, and the facets of this cone
correspond to $c_{IJ}^{K}=1$. P.~Belkale \cite{3} was first to
note that all inequalities (IJK) follows from those with
$c_{IJ}^{K}=1$, and in recent preprint A.~Knutson, T.~Tao, and
Ch.~Woodward \cite{23} proved their independence. In my original
paper \cite{19} condition  (2.3) appears in a weaker form
$$V_{N\gamma}\subset V_{N\alpha}\otimes V_{N\beta}\quad\text{for some }N>0,\tag{$2.3^\prime$}
 $$
and its equivalence to (2.3), known as {\it saturation
conjecture}, was later proved by A.~Knutson and T.~Tao \cite{22},
and in more general quiver context by H.~Derksen and J.~Weyman
\cite{6}.

Note that inequalities (IJK), although complete, are too numerous
to be practical for large $n$. That is why L-R rule, in its
different incarnations \cite{22, 11}, often provides a more
intuitive way to see possible spectra for sum of Hermitian
operators.

\subhead{Example~2.3}\endsubhead  Let $A$ be Hermitian matrix with integer
spectrum $\lambda(A):\,\,a_1\ge a_2\ge\ldots\ge a_n$ and $B\ge0$
be a nonnegative matrix of rank one with spectrum
$\lambda(B):\,\,b\ge0\ge\cdots\ge0$. Viewing the spectra as Young
diagrams, and applying L-R rule we find out that
$\lambda(A)\otimes\lambda(B)$ is a sum of diagrams
$\gamma:\,\,\,c_1\ge c_2\ge\cdots \ge c_n$ satisfying the
following {\it intrlacing inequalities}
$$c_1\ge a_1\ge c_2\ge a_2\cdots\ge c_n\ge a_n.$$
By Theorem 2.3 this implies {\it Cauchy interlacing theorem} for
spectra
$$\lambda_i(A)\le\lambda_i(A+B)\le\lambda_{i-1}(A),\quad
\text{rk}B=1,\quad B\ge 0,$$
known in mechanics as {\it
Rayleigh-Courant-Fisher} principle: Let mechanical system
$S^\prime$ is obtained from another one $S$, by imposing a linear
constraint, e.g. by fixing a point of a drum. Then spectrum of
$S$ separates spectrum of $S^\prime$.

\specialhead \noindent \boldLARGE 3. Toric bundles \endspecialhead

Historically Theorem 2.3 first appears as a byproduct of theory of
toric vector bundles and sheaves, originated in \cite{15, 17}. See
other expositions of the theory in \cite{21, 30}, and further
applications in \cite{16, 33}. Vector bundles form a cross point
at which the diverse subjects of this paper meet together.

\specialhead \noindent \boldlarge 3.1. Filtrations
\endspecialhead

To avoid technicalities let's consider the simplest case of
projective plane
$$
\Bbb P^2=\{(x^\alpha:x^\beta:x^\gamma)|x\in \Bbb C\}
$$
on which diagonal torus
$$
T=\{(t_\alpha : t_\beta : t_\gamma)| t\in \Bbb C^*\}\tag 3.1
$$
acts by the formula
$$
t\cdot x=(t_\alpha x^\alpha:t_\beta x^\beta:t_\gamma x^\gamma).
$$
Orbits of this action are vertices, sides and complement of the
coordinate triangle. In particular there is unique dense orbit,
consisting of points with nonzero coordinates.

The objects of our interest are $T$-{\it equivariant} (or {\it
toric} for short) vector bundles $\Cal E$ over  $\Bbb P^2$. This
means that $\Cal E$ is endowed with an action $T:\Cal E$ which is
linear on fibers and makes the following diagram commutative
$$ \CD \Cal E @
> t>> \Cal E \\  @V\pi VV @VV\pi V \\ \Bbb P^2 @>t >>\Bbb P^2
\endCD \;\;\qquad t\in T.$$
Let us fix a generic point $p_0\in \Bbb P^2$ not in a coordinate
line, and denote by
$$
E:=\Cal E(p_0)
$$
the corresponding generic fiber. There is no action of torus $T$
on the fiber $E$. Instead the equivariant structure produces some
distinguished subspaces in $E$ by the following construction. Let
us choose a generic point $p_\alpha\in X^\alpha $ in coordinate
line $ X^\alpha:x^\alpha=0$. Since  $T$-orbit of $p_0$ is dense in
$\Bbb P^2$, we can vary $t\in T$ so that $t p_0$ tends to
$p_\alpha$. Then for any vector  $e\in E=\Cal E(p_0)$, we have
$te\in \Cal E(tp_0)$ and can try the limit
$$
\lim_{tp_0\rightarrow p_\alpha}(te)
$$
which either exists or not.  Let us denote by $E^\alpha(0)$ the
set of vectors $e\in E$ for which the limit exists:
$$
E^\alpha(0):=\{e\in E|\lim_{tp_0\rightarrow
p_\alpha}(te)\;\;\text{ exists}\}.
$$
Evidently $E^\alpha(0)$ is a vector subspace of $E$, independent
of $p_0$ and $p_\alpha$.

An easy modification of the previous construction allows to define
for integer $m\in
\Bbb Z$, the subspace
$$
E^\alpha(m):=\left\{e\in E |\lim_{tp_0\rightarrow
p_\alpha}\left(\frac{t_\alpha}{t_\beta}\right)^{-m}\cdot (te)
\;\;\text{ exists}\right\}.
$$
Roughly speaking $E^\alpha(m)$ consists of vectors $e\in E$ for
which $te$ vanishes up to order $m$ as $tp_0$ tends to coordinate
line $X^\alpha$. The subspaces $E^\alpha(m)$ form a non-increasing
exhaustive $\Bbb Z$-filtration:
$$\aligned E^\alpha :\cdots
\supset E^\alpha(m-1)&\supset E^\alpha(m)\supset
E^\alpha(m+1)\supset\cdots,
\\E^\alpha(m)&=\,\,0,\text{ for $m\gg 0$},\\
E^\alpha(m)&=E,\text{ for $m\ll 0$}.\endaligned\tag3.2$$ Applying
this construction to other coordinate lines, we get a triple of
filtrations $E^\alpha$, $E^\beta$, $E^\gamma$ in generic fiber
$E=\Cal E(p_0)$, associated with  toric bundle $\Cal E$.

\proclaim{\enspace Theorem 3.1} The correspondence
$$
\Cal E\mapsto (E^\alpha,E^\beta,E^\gamma)\tag3.3
$$
establishes an equivalence between  category of toric vector
bundles on $\Bbb P^2$ and category of triply filtered vector
spaces.
\endproclaim
We'll use notation $\Cal E(E^\alpha,E^\beta,E^\gamma)$ for toric
bundle corresponding to triplet of filtrations
$E^\alpha,E^\beta,E^\gamma$.

\specialhead \noindent \boldlarge3.2. Stability \endspecialhead

The previous theorem tells that every property or invariant of a
vector bundle has its counterpart on the level of filtrations. For
application to spectral problems the notion of stability of a
vector bundle $\Cal E$ is crucial. Recall that $\Cal E\rightarrow
\Bbb P^2$ is said to be Mumford--Takemoto {\it stable\/} iff
$$\frac{c_1(\Cal F)}{\operatorname{rk}\Cal F}<
\frac{c_1(\Cal E)}{\operatorname{rk}\Cal E}\tag3.4$$
for every proper subsheaf $\Cal F\subset \Cal E$, and {\it
semistable} if weak inequalities hold. Here $c_1(\Cal
E)=\operatorname{deg}\det\Cal E$ is the first Chern class.
Donaldson theorem \cite{7} brings a deep geometrical meaning to
this seemingly artificial definition: Every stable bundle carries
unique Hermit-Einstein metric (with Ricci curvature proportional
to metric).
\proclaim{Theorem 3.2} Toric bundle $\Cal E=\Cal
E(E^\alpha,E^\beta,E^\gamma)$ is stable iff for every proper
subspace $F\subset E$ the following inequality holds
$$\frac1{\dim F}\sum\Sb \nu=\alpha,\beta,\gamma\\i\in\Bbb
Z\endSb i\dim F^{[\nu]}(i)<
\frac1{\dim E}\sum\Sb \nu=\alpha,\beta,\gamma\\i\in\Bbb
Z\endSb i\dim E^{[\nu]}(i)\tag3.5$$ where $F^{\nu}(i)=F\cap
E^\nu(i)$ is induces filtration with composition factors
$F^{[\nu]}(i)=F^\nu(i)/F^\nu(i+1)$.
\endproclaim

There is nothing surprising in this theorem since the sums in
(3.5) are just Chern classes of the corresponding toric bundles
and sheaves.
\subhead Remark~3.3\endsubhead Inequality (3.5) depends only on
{\it relative positions\/} of subspace $F\subset E$ with respect
to filtrations $E^\alpha, E^\beta,E^\gamma$, which are given by
three {\it Schubert cells} $s_\alpha,s_\beta,s_\gamma$. Hence we
have one inequality each time.
$$s_\alpha\cap s_\beta\cap s_\gamma\ne\emptyset.\tag3.6$$
For filtrations in general position (3.6) is equivalent to
nonvanishing of the product of {\it Schubert cycles\/}
${\sigma}_\alpha\cdot{\sigma}_\beta\cdot{\sigma}_\gamma\ne0$ in
cohomolgy ring of Grassmannian, and in this case  stability
inequalities (3.5) amount  to inequalities (IJK) of Theorem 2.1.

\specialhead \noindent \boldlarge 3.3. Back to Hermitian
operators \endspecialhead

Let now $E$ be Hermitian space and $H:E\rightarrow E$ be Hermitian
operator with {\it spectral filtration}
$$E^H(x)=\left(\vcenter{\hsize3.95cm\noindent sum of eigenspaces
of $H$ with eigenvalues at least $x$}\right).\tag3.7$$ The
operator can be recovered from the filtration using  {\it spectral
decomposition}
$$H=\int_{-\infty}^\infty xdP_H(x)$$
where $P_H(x)$ is orthogonal projector with kernel $E^H(x)$. So in
Hermitian space we have equivalence
$${\text{\it Hermitian operators }= \Bbb
R\text{\it-filtrations}}.$$ Let $H^\alpha$ be Hermitian operator
with spectral filtration $E^\alpha$. Its spectrum  depends only on
filtration $E^\alpha$, and we define
$\operatorname{Spec}E^\alpha:=\operatorname{Spec}H^\alpha$.
\proclaim{Theorem~3.3} Indecomposable triplet of $\Bbb
R$-filtrations $E^\alpha,E^\beta,E^\gamma$ is stable iff there
exists a Hermitian metric in $E$ such that the sum of the
corresponding Hermitian operators is a scalar
$$H^\alpha+H^\beta+H^\gamma=\text{\rm scalar}.\tag3.8$$
\endproclaim
This is a toric version of Donaldson theorem on existence of
Hermit--Einstein metric in stable bundles. Together with Theorem
3.2 it reduces solution of Hermitian spectral problem to stability
inequalities (3.5), which by remark 3.3 amounts to inequalities
(IJK) of Theorem 2.1.

See also Faltings talk \cite{9} on arithmetical applications of
stable filtrations.

\specialhead \noindent \boldlarge 3.4. Components of tensor
product \endspecialhead

In the previous section we explain that stability inequalities
(3.5) ($\Leftrightarrow$  (IJK)) via toric Donaldson-Yau theorem
solve Hermitian spectral problem. To relate this with tensor
product part of Theorem 2.1 we need another interpretation of the
stability inequalities via Geometric Invariant Theory \cite{26}.

Recall, that point $x\in \Bbb P(V)$ is said to be {\it GIT stable
} with respect to linear action $G:V$ if $G$-orbit of the
corresponding vector $\overline{x}\in V$ is closed and its
stabilizer is finite. Let
$$X=\Cal F^\alpha\times \Cal F^\beta\times \Cal F^\gamma$$
be product of three flag varieties of the same types as flags of
the filtrations $E^\alpha,E^\beta,E^\gamma$, and $\Cal L^\alpha$
be line bundle on the flag variety $\Cal F^\alpha$ induced by
character
$$\omega_\alpha:\text{diag}(x_1,x_2,\ldots,x_n)\mapsto
x_1^{a_1}x_1^{a_2}\cdots x_1^{a_n},$$ where $\alpha:\,\,a_1\ge
a_2\ge\cdots\ge a_n$ is the spectrum of filtration $E^\alpha$,
i.e. spectrum of the corresponding operator $H^\alpha$.
\proclaim{Observation 3.4} Vector bundle $\Cal E=\Cal
E(E^\alpha,E^\beta,E^\gamma)$ is stable iff the corresponding
triplet of flags
$$x=F^\alpha\times F^\beta \times F^\gamma\in
\Cal F^\alpha\times \Cal F^\beta\times \Cal F^\gamma=X\hookrightarrow
\Bbb P(\Gamma(X,\Cal L))$$
is a GIT stable point w.r. to group $\text{\rm SL}(E)$ and
polarization $\Cal L=\Cal L^\alpha\boxtimes\Cal
L^\beta\boxtimes\Cal L^\gamma$.
\endproclaim

This observation is essentially due to Mumford \cite{25}. Notice
that by Borel-Weil-Bott theorem \cite{5} the space of global
sections $\Gamma(\Cal F^\alpha,\Cal L^\alpha)=V_\alpha$ is just an
irreducible representation of $\operatorname{SL}(E)$ with highest
weight $\alpha$. Hence $\Gamma(X,\Cal L))=V_\alpha\otimes
V_\beta\otimes V_\gamma$. Every stable vector $\overline{x}$ can
be separated from zero by a $G$-invariant section of $\Cal L^N$.
Therefore triplet of flags in generic position is stable iff
$[V_{N\alpha}\otimes V_{N\beta}\otimes
V_{N\gamma}]^{\text{SL}(E)}\ne0$ for some $N\ge 1$. This proves
the last part of Theorem 2.1, modulo the saturation conjecture.

\specialhead \noindent \boldLARGE 4. Unitary operators and
parabolic bundles \endspecialhead

We have seen in the previous section that solution of the
Hermitian spectral problem amounts to stability condition for
toric bundles. A remarkable ramification of this idea was
discovered by S.~Angihotri and Ch.~Woodward \cite{2} for unitary
spectral problem.

Let $U\in \text{SU}(n)$ be unitary matrix with unitary spectrum
$$\varepsilon(U)=(e^{2\pi i\lambda_1},e^{2\pi i\lambda_2},\ldots,e^{2\pi
i\lambda_n}).$$ Let's normalize  exponents $\lambda_i$ as follows
$$\lambda(U):=\left\{\aligned
&\lambda_1\ge\lambda_2\ge\cdots\ge\lambda_n,\\
&\lambda_1+\lambda_2+\cdots+\lambda_n=0,\\
&\lambda_1-\lambda_n<1,\endaligned\right.\tag4.1$$ and, admitting
an abuse of language,  call $\lambda(U)$ {\it spectrum} of $U$.

\proclaim{Unitary spectral problem} {\rm Find possible spectra of
product $\lambda(UV)$, when spectra of the factors
$\lambda(U),\lambda(V)$ are given.}
\endproclaim

To state the result we need in {\it quantum cohomology}
$\text{H}^*_q(G_p^r)$ of Grassmannian $G_p^r$ of linear subspaces
of dimension $p$ and codimension $r$. This is an algebra over
polynomial ring $\Bbb C[q]$ generated by Schubert cycles
$\sigma_I$, $I\subset \{1,2,\ldots,n\},\,|I|=p$, $n=p+r$ with
multiplication given by the formula
$$\sigma_I*\sigma_J=\sum_{K,d}c_{IJ}^K(d)q^d\sigma_K$$
where structure constants $c_{IJ}^K(d)$ are defined as follows.
Let $G^p_q\hookrightarrow\Bbb P(\bigwedge^p\Bbb C^n)$ be Pl\"ucker
imbedding and
$$\varphi:\Bbb P^1\rightarrow G_p^r$$
be a rational curve of degree $d$ in Gressmanian $G^p_q\subset\Bbb
P(\bigwedge^p\Bbb C^n)$. One can check that $\varphi$ depends on
$\dim G^p_q+nd$ parameters. For fixed point $x\in\Bbb P^1$ the
condition $\varphi(x)\in \sigma_I$ imposes
$\text{codim\,}\sigma_I$ constraints on $\varphi$. Hence for
$$\text{codim\,}\sigma_I+\text{codim\,}\sigma_J+\text{codim\,}\sigma_K=\dim
G^r_p+nd$$ the numbers
$$(\sigma_I,\sigma_J,\sigma_K)_d=\#\{\varphi:\Bbb P^1\rightarrow G_p^r\mid \varphi(x)\in
\sigma_I,\varphi(x)\in
\sigma_J,\varphi(x)\in
\sigma_K,\deg\varphi=d\}$$
supposed to be finite. They are known as {\it Gromov -Witten
invariants} and related to the structure constants by the formula
$$c_{IJ}^K(d)=(\sigma_I,\sigma_J,\sigma_{K^*})_d$$
where $K^*=\{n+1-k\mid k\in K\}$. For $d=0$ they are just
conventional Littlewood--Richardson coefficients $c_{IJ}^K$.
\proclaim{Theorem 4.1} The following conditions are equivalent

{\rm i)} There exist unitary matrices $W=UV$ with given spectra
$\lambda(U),\lambda(V),\lambda(W)$.

{\rm ii)} The inequality

\hfill $\lambda_I(U)+\lambda_J(V)\le d+\lambda_K(W)${\hskip4cm
$\text{\rm (IJK)}_d$}

holds each time $c_{IJ}^K(d)\ne0$.
\endproclaim

\specialhead \noindent \boldlarge 4.1. Parabolic bundles
\endspecialhead As in the Hermitian case solution of the unitary
problem comes from its {\it holomorphic} interpretation in terms
of vector bundles. To explain the idea let's start with vector
bundle $\Cal E$ over {\it compact} Riemann surface $\overline{X}$
of genus $g\ge 2$. It has unique topological invariant $c_1(\Cal
E)=\deg\det\Cal E$, which for simplicity we suppose to be zero,
i.e. $\Cal E$ be topologically trivial. Narasimhan-Seshadri
theorem \cite{27} claims that every stable bundle carries unique
flat metric, and hence defines unitary monodromy representation
$$\rho_{\Cal E}:\pi_1(\overline{X},x_0)\rightarrow
\text{SU}(E),\quad E=\Cal E(x_0).$$
This gives rise to equivalence
$$\Cal M_g:=\left(\vcenter{\hsize2.3cm\noindent stable bundles
\newline of degree zero}\right)=\left(\vcenter{\hsize4cm\noindent irreducible
uitary representations $\rho:\pi_1\rightarrow\text{SU}(E)$
}\right).\tag4.2$$ This theorem is an ancestor of the
Donaldson-Yau generalization \cite{7} to higher dimensions, and
may be seen as a geometric version of Langlands correspondence.

In algebraic terms the theorem describes stable bundles in terms
of solution of equation
$$[U_1,V_1][U_2,V_2]\cdots[U_g,V_g]=1$$
in unitary matrices $U_i,V_j\in\text{SU}(E)$. This is not the
matrix problem we are currently interested in. To modify it let's
consider {\it punctured\/} Riemann surface
$X=\overline{X}\backslash\{p_1,p_2,\ldots,p_\ell\}$. It has
distinguished classes
$$\gamma_\alpha=\left(\text{small circle around $p_\alpha$}\right)$$
in fundamental group $\pi_1(X)$, and we can readily define an
analogue of  RHS of (4.2):
$$\Cal
M_g(\lambda^{(1)},\lambda^{(2)},\cdots,\lambda^{(\ell)})=\{\rho:\pi_1(X)\rightarrow
\text{SU}(E)\mid \lambda(\rho(\gamma_\alpha))=\lambda^{(\alpha)}\},\tag4.3$$
where $\lambda^{(\alpha)}$ is a given spectrum of monodromy around
puncture $p_\alpha$. C.~S.~Seshadri \cite{31} manages to find an
analogue of more subtle {\it holomorphic\/} LHS of (4.2) in terms
of so called {\it parabolic bundles\/}.

Parabolic bundle $\Cal E$ on $X$ is actually a bundle on
compactification $\overline{X}$ together with $\Bbb R$-filtration
in every special fiber $E^\alpha=\Cal E(p_\alpha)$ with support in
an interval of length $\le1$. The filtration is a substitution for
spectral decomposition of $\rho(\gamma_\alpha)$, cf. (4.1).
Seshadri also defines {\it (semi)stability\/} of parabolic bundle
$\Cal E$ by inequalities
$$\frac{\text{Par}\deg\Cal F}{\text{rk\,}{\Cal F}}\le
\frac{\text{Par}\deg\Cal E}{\text{rk\,}{\Cal E}},\quad\forall \Cal
F\subset\Cal E,\tag4.4$$ where the parabolic degree is given by
equation $\text{Par}\deg\Cal E=\deg\Cal E+\sum_{\alpha,i}
\lambda^{(\alpha)}_i$. Metha-Seshadri theorem \cite{24}
claims that every stable parabolic bundle $\Cal E$ on $X$ carries
unique flat metric with given spectra of monodromies
$\lambda(\gamma_\alpha)=\lambda^{(\alpha)}$. This gives a
holomorphic interpretation of the space (4.3)
$$\Cal
M_g(\lambda^{(1)},\lambda^{(2)},\cdots,\lambda^{(\ell)})=
\left(\vcenter{\hsize6cm\noindent
stable parabolic bundles of degree zero with given types of the
filtrations }\right).\tag4.5$$ In the simplest case of projective
line with three punctures (4.3) amounts to space of solutions of
equation $UVW=1$ in unitary matrices $U,V,W\in
\text{SU}(n)$ with given spectra. By Metha-Seshadry theorem
solvability of this equation is equivalent to stability
inequalities (4.4). In the case under consideration holomorphic
vector bundle $\Cal E$ on $\Bbb P^1$ is trivial, $\Cal
E=E\times\Bbb P^1$, and hence its subbundle $\Cal F\subset \Cal E$
of rank $p$ is nothing but a rational curve $\varphi:\Bbb
P^1\rightarrow G_p(E)$ in Grassmannian. This allows to write down
stability condition (4.4) in terms of  quantum cohomology, and
eventually arrive at Theorem 4.1.

\specialhead \noindent \boldLARGE 5. Further
ramifications\endspecialhead

The progress in Hermitian and unitary spectral problems open way
for solution of a variety of  others classical, and not so
classical, problems. Most of them, however, have no holomorphic
interpretation, and require different methods, borrowed from
harmonic analysis on homogeneous spaces, symplectic geometry, and
geometric invariant theory.

\specialhead \noindent \boldlarge 5.1. Multiplicative singular
value problem\endspecialhead

The problem in question is about possible singular spectrum
$\sigma(AB)$ of product of complex matrices with given singular
spectra $\sigma(A)$ and $\sigma(B)$. Recall, that singular
spectrum of complex matrix $A$ is spectrum of its radial part
$\sigma(A):=\lambda(\sqrt{A^*A})$.

For a long time it was observed that every inequality for
Hermitian problem has a {\it multiplicative} counterpart for the
singular one. For example multiplicative version of Weyl's
inequality $\lambda_{i+j-1}(A+B)\le\lambda_i(A)+\lambda_j(B)$ is
$\sigma_{i+j-1}(AB)\le\sigma_i(A)\sigma_j(B)$. The equivalence
between these two problems was conjectured by R.~C.~Thompson, and
first proved by the author \cite{20} using harmonic analysis on
symmetric spaces. Later on A.~Alekseev, E.~Menreken, and
Ch.~Woodward \cite{1} gave an elegant conceptual solution based on
Drinfeld's Poisson-Lie groups \cite{8}. Here is a precise
statement for classical groups.
\proclaim{Theorem\enspace 5.1}
 Let $\Bbb G$ be one of the classical groups $\text{\rm SL}(n,\Bbb
C)$, $\text{\rm SO}(n,\Bbb C)$, or $\text{\rm Sp}(2n,\Bbb C)$ and
$\Bbb L$ be the corresponding compact Lie algebra of traceless
skew Hermitian complex, real, or quaternionic $n\times n$ matrices
respectively. Then the following conditions are equivalent
\roster
\item There exist $A_i\in \Bbb G$ with given singular spectra
$\sigma(A_i)=\sigma_i$ and
$$A_1A_2\cdots A_N=1.$$
\item There exist $H_i\in\Bbb L$ with spectra
$\lambda(H_i)=\sqrt{-1}\log\sigma_i$ and
$$H_1+H_2+\cdots+H_N=0.$$
\endroster\endproclaim

Note, however, that neither of the above approaches solve the
singular problem per se, but reduces it to Hermitian one. Both of
them suggest that all three problems must be treated in one
package. More precisely, every compact simply connected group $G$
give birth to three symmetric spaces
\roster
\item"$\bullet$" The group  $G$ itself,
\item"$\bullet$" Its Lie algebra $L_G$,
\item"$\bullet$" The dual symmetric space $H_G=G^{\Bbb C}/G$,
\endroster
of positive, zero, and negative curvature, and to three ``spectral
problems" concerned with support of convolution of $G$ orbits in
these spaces, see \cite{20} for details. For $G=\text{SU}(n)$ we
return to the package of unitary, Hermitian, and singular
problems.

The first two problems may be effectively treated in framework of
vector bundles with structure group $G$, as explained in sections
2--4. Many flat, i.e. additive ``spectral problem" has been solved
by A.~Berenstein and R.~Sjammar in a very general setting
\cite{4}.

\specialhead \noindent \boldlarge 5.2. Other symmetric
spaces\endspecialhead

As an example of  unresolved problem let's consider symmetric
spaces associated with different incarnations of Grassmannian
\roster
\item"$\bullet$" Compact $U(p+q)/U(p)\times U(q)$,
\item"$\bullet$" Flat $\text{Mat}(p,q)$ = coomplex $p\times q$ matrices,
\item"$\bullet$" Hyperbolic $U(p,q)/U(p)\times U(q)$.
\endroster
In compact case the corresponding spectral problem is about
possible angles between three $p$-subspaces $U, V, W\subset\Bbb
H^n$ in Hermitian space $\Bbb H^n$ of dimension $n=p+q$, $p\le q$.
The {\it Jordan angles}
$$\widehat{UV}=(\varphi_1,\varphi_2,\ldots,\varphi_p),\quad\quad 0\le \varphi\le \frac\pi2$$
between subspaces $U,V$ are defined via spectrum of product of
orthogonal projectors $\pi_{UV}:U\rightarrow V$ and
$\pi_{VU}:V\rightarrow U$
$$\lambda(\sqrt{\pi_{UV}\pi_{VU}}):
\cos\varphi_1\ge\cos\varphi_2\ge \cdots\ge\cos\varphi_p\ge0.$$
Yu.~Neretin \cite{28} proved Lidskii type inequalities
\footnote{He actually deals with {\it real} Grassamnnian.}
for angles $\widehat{UV},\widehat{VW},\widehat{WU}$, and
conjectured that other inequalities are the same as in the
Hermitian case. Note, however, that the unitary triplet suggests
existence of nonhomogeneous ``quantum" inequalities, e.g. sum of
angles of a spherical triangle is $\le\pi$.

In flat case the problem is about relation between singular
spectra of $p\times q$ matrices $\sigma(A-B)$, $\sigma(B-C)$,
$\sigma(C-A)$. This {\it additive} singular problem was resolved
by O'Shea and Sjamaar \cite{29}.

In hyperbolic case the question is about angles between maximal
{\it positive\/} subspaces $U,V,W\subset \Bbb H^{pq}$ in Hermitian
space of signature $(p,q)$. They are defined by equation
$$\lambda(\sqrt{\pi_{UV}\pi_{VU}}):
\cosh\varphi_1\ge\cosh\varphi_2\ge \cdots\ge\cosh\varphi_p\ge1.$$
Again our experience with the unitary triplet suggests that the
exponential map establishes a Thompson's type correspondence
between O'Shea-Sjamaar inequalities for additive singular problem
and that of for hyperbolic angles.

\specialhead \noindent \boldlarge 5.3. P-adic spectral
problems\endspecialhead

There is also a nonarchimedian counterpart of this theory, which
deals with classical Chevalley groups $\Bbb G_p=\text{\rm
SL}(n,\Bbb Q_p)$, $\text{\rm SO}(n,\Bbb Q_p)$, or $\text{\rm
Sp}(2n,\Bbb Q_p)$ over $p$-adic field $\Bbb Q_p$ and their maximal
compact subgroups $\Bbb K_p=\text{\rm SL}(n,\Bbb Z_p)$, $\text{\rm
SO}(n,\Bbb Z_p)$, or $\text{\rm Sp}(2n,\Bbb Z_p)$ respectively.
Double coset $\Bbb K_pg\Bbb K_p$ may be treated as a complete
invariant of lattice $L=gL_0$, $L_0=\Bbb Z_p^{\oplus n}$ with
respect to $\Bbb K_p$. We call lattice $L=gL_0$ {\it unimodular},
{\it orthogonal} or {\it symplectic} if respectively
$g\in\text{SL}(n,\Bbb Q_p)$, $g\in\text{SO}(n,\Bbb Q_p)$ or
$g\in\text{Sp}(2n,\Bbb Q_p)$.

It is commonly known that in the unimodular case there exists a
basis $e_i$ of $L_0$ such that $\tilde e_i=p^{a_i}e_i$ form a
basis of $L$ for some $a_i\in \Bbb Z$. We define {\it index}
$(L:L_0)$ by
$$(L:L_0)=(p^{a_1},p^{a_2},\ldots,p^{a_n}),\quad\quad
a_1\ge a_2\ge\cdots\ge a_n.\tag5.1$$ Notice that unimodularity
$g\in\text{SL}(n,\Bbb Q_p)$ implies $a_1+a_2+\cdots+a_n=0$.

The index $(L:L_0)$ of an orthogonal or a symplectic lattices has
extra symmetries. In orthogonal case  we may choose the above
basis $e_i$ of $L_0$ to be {\it neutral}, in which case the
quadratic form becomes
$$\sum_{1-n}^{n-1}x_ix_{-i},\quad i\equiv n-1 \mod 2.$$
Then the index takes the form
$$(L:L_0)=(p^{a_{n-1}},p^{a_{n-3}},\ldots,p^{a_{3-n}},p^{a_{1-n}}),\tag5.2
$$ where $a_{n-1}\ge a_{n-3}\ge\ldots\ge a_{3-n}\ge a_{1-n}$, and
$a_{-i}=-a_i$.

Similarly,  for {\it symplectic} lattice $L$ we can choose
symplectic basis $e_i,f_j$ of $L_0$ such that $\tilde
e_i=p^{a_i}e_i$ and $\tilde f_j=p^{-a_j}f_j$ form a basis of $L$.
In this case we have
$$(L:L_0)=(p^{a_n},p^{a_{n-1}},\ldots,p^{a_1},p^{-a_1},
\ldots,p^{-a_{n-1}},p^{-a_n}),\tag5.3$$
with $a_n\ge a_{n-1}\ge,\ldots,\ge a_1\ge 0$.

Notice that the spectra (5.1)-(5.3) have the same symmetry, as
singular spectrum $\sigma(A)$ of a matrix $A\in\Bbb G$ in the
corresponding classical {\it complex} group.
\proclaim{Theorem\enspace 5.2}
The following conditions are equivalent
\roster\item
There exists a sequence of (unimodular, orthogonal, symplectic)
lattices $$L_0,L_1,\ldots,L_{N-1},L_N=L_0$$ of given indices
$\sigma_i=(L_i:L_{i-1})$.
\item
The indices $\sigma_i$ satisfy the equivalent conditions of
Theorem 5.1 for the corresponding complex group $\Bbb G$.
\endroster\endproclaim
We'll give proof elsewhere. The theorem is known for the
unimodular lattices, see \cite{10}.

\specialhead \noindent \boldlarge 5.4. Final
remarks\endspecialhead

In the talk I try to trace the flaw of ideas from the theory of
vector bundles to spectral problems. It seems C.~Simpson \cite{32}
was the first to note that vector bundles technic has nontrivial
implications in linear  algebra. He proved that product
$C_1C_2\cdots C_N$ of conjugacy classes $C_i\subset
\text{SL}(n,\Bbb C)$ is dense in $\text{SL}(n,\Bbb C)$ iff
$$\aligned\dim C_1+\dim C_2+\cdots+\dim C_N&\ge (n+1)(n-2),
\\r_1+r_2+\cdots+r_N&\ge n,
\endaligned\tag5.4$$
where $r_i$ is maximal codimension of root space of a matrix
$A_i\in C_i$. This problem was suggested  by P.~Deligne, who noted
that under condition
$$\dim C_1+\dim C_2+\cdots+\dim C_N=2n^2-2$$
an irreducible solution of equation $A_1A_2\cdots A_N=1$, $a_i\in
C_i$ is unique up to conjugacy, see book of N.~Katz \cite{14} on
this rigidity phenomenon.

I think that inverse applications to moduli spaces of vector
bundles are sill ahead. One may consider {\it polygon spaces}
\cite{18, 12} as a toy example of this feedback, corresponding to
toric $2$-bundles. A similar space of {\it spherical polygons} in
$\Bbb S^3$ with given sides is a model for moduli space of flat
connections in punctured Riemann sphere. Its description is a
challenge problem.

There are many interesting results, e.g. infinite dimensional
spectral problems, which fall out of this survey. I refer to
Fulton's paper \cite{10} for missing details.

\specialhead \noindent \boldLARGE References \endspecialhead
\widestnumber\key{[XX]}

\ref\key1 \by A.~Alekseev, E.~Meinrenken, \& C.~Woodward \paper
{\rm Linearization of Poisson actions and singular values of
matrix product} \jour {\it Ann. Inst. Fourier (Grenoble),} \vol
{\rm 51} \issue6 \yr 2001 \pages 1691--1717
\endref

\ref\key2 \by S.~Angihotri \& C.~Woodward \paper {\rm Eigenvalues
of products of unitary matrices and quantum Schubert calculus}
\jour {\it Math. Res. Letters,} \vol {\rm 5} \yr 1998 \pages
817--836
\endref

\ref\key3 \by P.~Belkale \paper{\rm Local systems on $\Bbb {P}\sp
1-{S}$ for ${S}$ a finite set} \jour{\it Compositio Math.,}
\yr{2001} \issue 1 \vol{\rm 129} \pages{67--86}
\endref

\ref\key4 \by A.~Berenstein \& R.~Sjamaar \paper{\rm Coadjoint
orbits, moment polytopes, and the Hilbert--Mumford criterion}
\jour{\it J. Amer. Math. Soc.,} \vol{\rm 13} \yr 2000 \issue2
\pages 433--466
\endref

\ref \key 5 \by R.~Bott \paper {\rm Homogeneous vector bundles}
\jour {\it Ann. of Math.,} \vol{\rm 66} \yr 1957 \pages 203--248
\endref

\ref\key 6 \by{H.~Derksen \& J.~Weyman} \paper{\rm
Semi-invariants of quivers and saturation for
Littlewood-Richardson theorem} \jour{\it J. Amer. Math. Soc.,}
\vol{\rm 13} \issue3 \yr2000 \pages{467--479}
\endref

\ref \key 7 \by S.~K.~Donaldson \paper {\rm Infinite
determinants, stable bundles and curvature} \jour {\it Duke Math.
J.,} \vol {\rm 54} \yr 1987 \pages 231--247
\endref

\ref\key 8 \by V.~G.~Drinfeld \paper{\rm Quantum groups} \inbook
{\it Proceedings of the International Congress of Mathematicians}
\vol 1,2 (Berkeley, 1986) \publ Amer. Math. Soc. \publaddr
Providence, RI \yr 1987, 798--820
\endref

\ref\key 9 \by G.~Faltings \paper{\rm Mumford-Stabilit\"at in der
algebraischen Geometrie} \inbook{\it Proceedings of the
International Congress of Mathematiciens} \vol{\rm 1,2,
(Z\"urich, 1994)} \publ Birk\"hauser \publaddr Basel \yr 1995,
648--655
\endref

\ref\key{10} \by W.~Fulton \paper{\rm Eigenvalues, invariant
factors, highest weights, and Schubert calculus} \jour {\it Bull.
Amer. Math. Soc.,} \vol{\rm 37} \issue3 \yr 2000 \pages 209--249
\endref

\ref\key 11
\by O.~Gleizer \& A.~Postnikov
\paper{\rm Littlewood-{R}ichardson coefficients via {Y}ang-{B}axter
              equation}
\jour{\it Internat. Math. Res. Notices}
\issue14
\pages741--774
\yr 2000
\endref

\ref\key 12 \by {J.-C.~Hausmann \& A.~Knutson} \paper{\rm The
cohomology ring of polygon spaces} \jour{\it Ann. Inst. Fourier
(Grenoble),} \vol{\rm 48} \yr{1998} \issue{1} \pages{281--321}
\endref

\ref\key 13 \by A.~Horn \paper{\rm Eigenvalues of sum of Hemitian
matrices} \jour{\it Pacific J. Math.,} \vol{\rm 12} \yr 1962
\pages225--241
\endref

\ref\key 14
\by N.~M.~Katz
\book{\it Rigid local systems}
\publ Princeton University Press
\publaddr Princeton
\yr 1996
\endref

\ref\key 15 \by A.~A.~Klyachko \paper{\rm Equivariant bundles on
toric varieties} \jour{\it Izv. Akad. Nauk SSSR Ser. Mat.,}
\vol{\rm 53} \issue 5 \pages1001--1039 \yr 1989 \lang Russian
\transl \jour{\it Math. USSR-Izv.,} \vol{\rm 35} \issue 2 \yr1990
\pages63--64
\endref

\ref\key{16} \by A.~A.~Klyachko \paper {\rm Moduli of vector
bundles and class numbers} \jour {\it Functional. i Priozhen.,}
\vol{\rm 25} \yr 1991 \pages 81--83 \lang Russian \transl \jour
{\it Funct. Anal. Appl.,} \vol{\rm25} \yr 1991 \issue 1 \pages
67--69
\endref

\ref\key{17}
\by A.~A.~Klyachko
\paper{\rm Vector bundles and torsion free
sheaves on the projective   plane}
\jour {\it Preprint Max-Planck-Institute fur
Mathematik}
\vol {\rm MPI/91-59,}
\yr 1991
\endref

\ref\key 18 \by A.~A.~Klyachko \paper {\rm Spatial polygons and
stable configurations of points in the projective line}
\inbook{\it Algebraic geometry and its applications (Yaroslavl,
1992)}  \publ Vieweg \publaddr Braunschweig \yr 1994, 67--84
\endref

\ref\key{19} \by A.~A.~Kyachko \paper{\rm Stable bundles,
repesentation theory and Hermitian operators} \jour {\it Selecta
Mathematica,} \vol{\rm 4} \yr 1998 \pages 419--445
\endref

\ref\key{20} \by A.~A.~Klyachko \paper {\rm Random walks on
symmetric spaces and and inequalities for matrix spectra} \jour
{\it Linear Algebra Appl.,} \vol{\rm 319} \yr 2000 \issue 2-3
\pages 37--59
\endref

\ref\key 21 \by A.~Knutson \& E.~Sharp \paper{\rm Sheaves on
toric varieties for physics} \jour{\it Adv. Theor. Math. Phys.,}
\vol{\rm 2} \yr1998 \issue 4 \pages873--961
\endref

\ref\key 22
\by{A.~Knutson \& T.~Tao}
\paper{\rm The honeycomb model of $\text{{G}{L}}(n,{\Bbb {C}})$ tensor
              products. {I}. {P}roof of the saturation conjecture}
\jour{\it J. Amer. Math. Soc.,} \vol{\rm 12} \yr1999 \issue2
\pages1055--1090
\endref

\ref\key 23
\by {A.~Knutson, T.~Tao \& Ch.~Woodward}
\paper{\rm The honecomb model for {$\text{\rm GL}(n,\Bbb{C})$}
 tensor products {II}: Facets of {Littlewood}--{Richardson} cone}
\jour{\it Preprint}
%month={february},
\yr{2001}
\endref

\ref \key 24 \by V.~B.~Metha \& C.~S.~Seshadri \paper {\rm Moduli
of vector bundles on curves with parabolic sructure} \jour{\it
Math. Ann.,} \vol{\rm  258} \yr 1980 \pages 205-239
\endref

\ref \key 25 \by D.~Mumford \paper {\rm Projective invariants of
projective structures} \inbook Proc. Int. Congress of Math.
Sockholm, 1963 \publ Almquist \& Wiksells \publaddr Uppsala \yr
1963, 526--530
\endref

\ref \key 26
\by D.~Mumford, J.~Fogarty, \& F.~Kirwan
\book Geometric invariant theory
\publ Springer
\publaddr Berlin
\yr 1994
\endref

\ref \key 27 \by M.~S.~Narasimhan \& C.~S.~Seshadri \paper {\rm
Stable and unitary vector bundles on a compact Riemann surface}
\jour {\it Ann. Math.,} \vol {\rm 82} \yr 1965 \pages 540--567
\endref

\ref\key{28} \by Yu.~Neretin \paper {\rm On Jordan angles and
triangle inequality in Grassmannian} \jour {\it Geom. Dedicata,}
\vol{\rm 86} \issue{1-3} \yr 2001 \pages 81--92
\endref

\ref\key{29} \by L.~O'Shea \& R.~Sjamaar \paper {\rm Moments maps
and Riemannian symmetric pairs} \jour {\it Math. Ann.,} \vol{\rm
317} \issue3 \yr 2000 \pages 415--457
\endref

\ref\key 30 \by M.~Perling \paper{\rm Graded rings and
equivariant sheaves on toric varieties} \jour{\it Preprint Univ.
Kaiserslauten,} \yr 2001
\endref

\ref \key 31 \by C.~S.~Seshadri \paper{\rm Moduli of vectir
bundles on curves with parabolic structures} \jour{\it Bull.
Amer. Math. Soc.,} \vol {\rm } \yr 1977 \pages 124--126
\endref

\ref\key 32 \by C.~T. Sympson \paper {\rm Product of Matrices}
\inbook {\it Differential geometry, global analysis, and
topology}, Canadian Math. Soc. Conf. Proc. \vol{\rm 12} \yr 1992,
157--185 \publ AMS \publaddr Providence RI
\endref

\ref\key 33 \by C.~Vafa \& E.~Witten \paper{A strong coupling
test of ${S}$-duality} \jour{\it Nuclear Phys. B,} \vol{\rm 431}
\issue{1-2} \yr{1994} \pages{3--77}
\endref

\ref \key 34 \by H.~Weyl \paper {\rm Das asymptotischer
Verteilungsgesetz der Eigenwerte lineare partialler
Differentialgleichungen} \jour {\it Math. Ann.,} \vol {\rm 71}
\yr 1912 \pages 441--479
\endref

\ref \key 35 \by H.~Weyl \paper {\rm Ramifications, old and new,
of the eigenvalue problem} \jour {\it Bull. Amer. Math. Soc.,}
\vol {\rm 56} \yr 1950 \pages 115--139
\endref

\enddocument